%
\magnification=\magstep1
\baselineskip = 20pt
\raggedbottom
\raggedright
\nopagenumbers
\def\BX{$\diamondsuit$}

\font\smallrm=cmr8
\font\smcaps=cmcsc10

\font\smcaps=cmcsc10
\font\ti=cmss12


%


\def\var{\hbox{\rm var}}
\def\cov{\hbox{\rm cov}}

%
%

\newcount\eqnumber
\def\eqdef#1{\global\advance\eqnumber by 1
        \expandafter\xdef
                \csname#1eqref\endcsname{\the\eqnumber}%
        \immediate\write\reffile{\string\def
                \expandafter\string\csname#1eqref\endcsname
                        {\the\eqnumber}}%
        \eqno
        \eqprint{\the\eqnumber}}%

\def\eqref#1{%
        \expandafter \ifx \csname#1eqref\endcsname \relax
                \ifxrefwarning \message{Undefined equation label '#1'.}
                \expandafter\def\csname#1eqref\endcsname{00}%
        \else \eqprint{\csname#1eqref\endcsname}%
        \fi}%
\def\eqprint#1{(#1)}%


\newwrite\reffile \newif\ifreffileopened
\def\openreffile{\ifreffileopened\else
        \reffileopenedtrue
        \immediate\openout\reffile = \jobname.aux
        \fi}%
\def\readreffile{%
        \testfileexistence{aux}%
        \iffileexists
                \begingroup
                        \@setletters
                        \input \jobname.aux
                \endgroup
        \else
                \message{No cross-reference file; I won't give you
                        warnings about undefined labels.}%
                \xrefwarningfalse
        \fi
        \openreffile}%
\def\@setletters{%
        \catcode`_ = \letter \catcode`+ = \letter
        \catcode`- = \letter \catcode`@ = \letter
        \catcode`0 = \letter \catcode`1 = \letter
        \catcode`2 = \letter \catcode`3 = \letter
        \catcode`4 = \letter \catcode`5 = \letter
        \catcode`6 = \letter \catcode`7 = \letter
        \catcode`8 = \letter \catcode`9 = \letter
        \catcode`( = \letter \catcode`) = \letter}%

%

\def\bib{\noindent\hangindent=15pt}
\def\section#1{\centerline{\ti #1\hfil}}

\tolerance=10000000

\def\THEOREMLM{ 1}
\def\THEOREMFDWN{ 2}
\def\THEOREMSDAP{ 3}
\def\THEOREMSS{ 4}
\def\THEOREMDUAL{ 5}
\def\THEOREMHYPER{ 6}
\def\THEOREMGLM{ 7}
\def\THEOREMGLMB{ 8}
\def\THEOREMGAP{ 9}

\pageno=1
\centerline{\bf HYPERBOLIC DECAY TIME SERIES}
\vskip 1in
\centerline{A.I. McLeod}
\vskip 1in
\centerline{Department of Statistical \& Actuarial Sciences}
\centerline{The University of Western Ontario}
\centerline{London, Ontario  N6A 5B7}
\centerline{Canada}
\vskip 1in
\centerline{\it October 1997}
\vskip 0.5in
\centerline{McLeod, A.I. (1998), Hyperbolic decay time series,}
\smallskip
\centerline{{\it The Journal of Time Series Analysis} 19, 473-484.}
\vfill\eject
\centerline{\bf Abstract}

\headline={\ifodd\pageno\rightheadline \else\leftheadline\fi}
\def\rightheadline{\smallrm\hfil\folio}
\def\leftheadline{\smallrm\folio\hfil}

Hyperbolic decay time series such as,
fractional Gaussian noise (FGN) or
fractional autoregressive moving-average (FARMA)
process, each exhibit two distinct types of behaviour:
strong persistence or antipersistence.
Beran (1994) characterized the family of strongly persistent
time series.
A more general family of hyperbolic decay time series is
introduced and its basic properties are characterized
in terms of the autocovariance and spectral density functions.
The random shock and inverted form representations
are derived.
It is shown that every strongly persistent series
is the dual of an antipersistent series and vice versa.
The asymptotic generalized variance of hyperbolic decay
time series with unit innovation variance is shown to be infinite
which implies that the variance of the
minimum mean-square error one-step linear predictor
using the last $k$ observations decays slowly to the innovation
variance as $k$ gets large.


\vfill\eject

\section{1. INTRODUCTION}
Let $Z_t, \, t=1,2,\ldots $ denote a covariance stationary,
purely nondeterministic
time series with mean zero and with autocovariance function,
$\gamma_Z(k) = \cov(Z_t, Z_{t-k})$.
As is discussed by Beran (1994),    many long memory processes
such as the FGN (Mandelbrot, 1983) and
FARMA (Granger and Joyeux, 1980; Hosking, 1981)
may be characterized by the property that
$k^{\alpha} \gamma_Z(k) \to  c_{\gamma}$ as $k \to \infty$,
for some  $\alpha \in (0,1)$ and $c_{\gamma} > 0$.
Equivalently,
$$
\gamma_Z(k) \sim  c_{\gamma}\, k^{-\alpha}.
\eqdef{LongMemory}
$$
As noted in Box and Jenkins (1976),
the usual stationary ARMA models
on the other hand are exponentially damped since
$\gamma_Z(k) = O(r^k)$, $r \in (0,1)$.

Beran (1994, p.42) shows that   an equivalent characterization
of strongly persistent time series is
$$f_Z(\lambda) \sim c_f \lambda^{\alpha-1} \quad {\rm as}\quad
 \lambda \rightarrow 0,
\eqdef{SDFOrigin}
$$
where
$\alpha \in (0,1)$, $c_f > 0$ and
$f_Z(\lambda)$ is the spectral density function given by
$f_Z(\lambda) =  \sum \gamma_Z(k) e^{-i k \lambda} /({2\pi})$.
Theorem \THEOREMLM\  below summarizes some results
stated without proof in Beran (1994, Lemma 5.1).
Since not all time series satisfying eq. \eqref{LongMemory} or
\eqref{SDFOrigin} are invertible, the restriction to invertible
processes is required.

{\smcaps Theorem \THEOREMLM .\/}
{\it The time series $Z_t$ satisfying \eqref{LongMemory}
or \eqref{SDFOrigin}
may be written in random shock form as
$Z_t = A_t + \sum \psi_{\ell} A_{t-\ell} $,
where
$\psi_{\ell} \sim   c_{\psi}\, \ell^{-(1+\alpha)/2}$,
$c_{\psi} > 0$ and $A_t$ is white noise.
Assuming that $Z_t$ is invertible, the inverted form may be written,
$Z_t = A_t + \sum\limits \pi_{\ell} Z_{t-\ell}$,
where
$\pi_{\ell} \sim    c_{\pi}\, \ell^{-(3-\alpha)/2}$,
$c_{\pi} > 0$ and $A_t$ is white noise.
\/}

{\smcaps Proof.\/}
By the Wold Decomposition, any purely nondeterministic time series may
be written in random shock form.
Now assume the random shock coefficients specified in the theorem and
we will derive \eqref{LongMemory}.
Assuming $\var(A_t) = 1$,
$\gamma_Z(k) = \psi_k + \sum \psi_h \psi_{h+k}$
$$\eqalign{\gamma_Z(k) &\sim
 \psi_k +  c_{\psi}^2 \sum_{h=1}^\infty h^{-(1+\alpha)/2} (h+k)^{-(1+\alpha)/2} \cr
&\sim
 \psi_k +  c_{\psi}^2 \int_{1}^\infty h^{-(1+\alpha)/2} (h+k)^{-(1+\alpha)/2} dh
 + R_k \cr
}$$
where the last step used
the Euler summation formula (Graham, Knuth and Patashnik, 1989, 9.78, 9.80)
and
$$ R_k = \{-{1\over 2} F(h)
+ {1\over 12} F^{\prime}(h)
+ {\theta\over 720} F^{\prime\prime\prime}(h) \} \Bigl|^\infty_1,
$$
where $\theta \in (0,1)$ and
$F(h) = h^{-(1+\alpha)/2} (h+k)^{-(1+\alpha)/2}.$
It is easily shown that $k^\alpha R_k \to 0$ as $k \to \infty$.
Hence,
$$\eqalign{\gamma_Z(k) &\sim
 \psi_k +  c_{\psi}^2 \int\limits_{1}^\infty
   h^{-(1+\alpha)/2} (h+k)^{-(1+\alpha)/2} dh \cr
&\sim
 \psi_k +  k^{-\alpha} c_{\psi}^2 \int\limits_{1/k}^\infty
   x^{-\beta} (x+1)^{-\beta} dx, \cr
}$$
where $\beta = (1+\alpha)/2$.
Using {\it Mathematica\/},
$$ \int\limits_{0}^\infty x^{-\beta} (x+1)^{-\beta} dx =
{{{2^{2\,\beta}}\,{\Gamma}(1 - \beta)\,{\Gamma}(-{1\over 2} + \beta)}\over
   {4\,{\sqrt{\pi }}}}
,$$
so \eqref{LongMemory} now follows with
$c_{\gamma} = c_{\psi}^2
{{{2^{\alpha-1}}\,{\Gamma}((1 - \alpha)/2)\,
{\Gamma}(\alpha/2)}/
   {\sqrt{\pi }}}
$, where $\Gamma(\bullet)$ is the gamma function.
This shows that $\psi_k$ is a possible factorization of $\gamma_k$
and that sufficies to establish that
$Z_t = A_t + \sum \psi_{\ell} A_{t-\ell} $.

For any stationary invertible linear process, $Z_t$,
$$\gamma_Z(k) = \sum\limits_{h=1}^\infty \pi_h \gamma_Z(k-h).
\eqdef{LinearProcess}
$$
Assume $\gamma_Z(k)$ satisfies eq. \eqref{LongMemory} and that
$\pi_{\ell} \sim    c_{\pi} \ell^{-(3-\alpha)/2}$
then we will show that eq. \eqref{LinearProcess} is satisfied.
$$\eqalign{\gamma_Z(k)
&\sim
\gamma_Z(0) \pi_k
+c \sum_{h=1}^{k-1} h^{-3/2+\alpha/2} (k-h)^{-\alpha}
+c \sum_{h=k+1}^{\infty} h^{-3/2+\alpha/2} (h+k)^{-\alpha},\cr}$$
where $c = c_{\pi} c_{\gamma}$.
Now
$\gamma_Z(0) \pi_k/\gamma_Z(k) \sim 0$ so the first term
will drop out.
In the second term, for $k>>h$, $(k-h)^{-\alpha} \sim k^{-\alpha}$ and
$$\sum_{h=1}^{k-1} h^{-3/2+\alpha/2} \sim H_{\beta}\quad
{\rm as}\ \ k \to \infty,$$
where
$H_{\beta} = \sum\limits_{h=1}^\infty h^{-\beta} < \infty$,
$\beta = 3/2-\alpha/2$.
In the final term, when $h>>k$, $(h+k)^{-\alpha} \sim h^{-\alpha}$,
so
$$
\eqalign{\sum_{h=k+1}^{\infty} h^{-3/2+\alpha/2} (h+k)^{-\alpha}
&\sim
 \sum\limits_{h=k+1}^{\infty} h^{-3/2 - \alpha/2} \cr
&\sim
 \int\limits_{k+1}^{\infty} h^{-3/2 - \alpha/2} dh \cr
&\sim
 (k+1)^{-(1 + \alpha)/2}. \cr
}
$$
Again the last step uses the Euler Summation Formula.
Thus the final term is smaller asymptotically smaller than $\gamma_k$.
This establishes the asymptotic equivalence of the left-hand side
and the right-hand side of eq. \eqref{LinearProcess}
and the theorem since $\gamma_Z(k)$ uniquely determines
the coefficients $\pi_{\ell}$ in the inverted model. \BX

The FARMA model of order $(p,q)$
(Granger and Joyeux, 1980; Hosking, 1981)
may be defined by the equation,
$$
\phi(B) (1-B)^d Z_t = \theta(B) A_t,
\eqdef{FARMA}
$$
where
$|d|<0.5$,
$A_t$ is white noise with variance $\sigma_A^2$,
$\phi(B) = 1-\phi_1 B -\ldots -\phi_p B^p$,
and
$\theta(B) = 1-\theta_1 B -\ldots -\theta_q B^q$.
For stationarity and invertibility it is assumed that
all roots of $\phi(B) \theta(B) = 0$ are outside the
unit circle and $|d|<0.5$.
The series is strongly persistent or antipersistent according
as $0<d<0.5$ or $-0.5<d<0$.
The special case where $p=q=0$ is known as fractionally
differenced white noise.

Antipersistent series may arise in practice when modelling nonstationary
time series.
As suggested by Box \& Jenkins (1976) a nonstationary time series
can often be made stationary by differencing the series until
stationarity is reached.
Sometimes the resulting stationary time series may be usefully
modelled by an antipersistent form of the FARMA model.
An illustrative example is provided by the annual U.S. electricity
consumption data for 1920--1970.
Hipel and McLeod (1994, pp.154--159) modelled the square-root consumption
using an ARIMA(0,2,1) but a better fit is obtained by
modelling the second differences of the square-root consumption
as fractionally differenced white noise with $d=-0.4477 \pm 0.1522$ sd.
The AIC for the latter model is 1011.5 as compared with 1020.4.
Diebold and Rudebusch (1989)
and
Beran (1995)
also used this approach for modelling nonstationary data.

The determinant of the covariance matrix of $n$ successive observations,
$Z_t,\ t=1,\ldots,n$, is denoted by
$G_Z(n) = \det(\gamma_Z(i-j))$.
It will now be shown in Theorem \THEOREMFDWN\
that for fractionally differenced
white noise,
$g_Z(n) = \sigma_A^{-2n} G_Z(n) \rightarrow \infty$
as
$n \rightarrow \infty$, where
$0<\sigma_A^2<\infty$, is the innovation variance given by
Kolmogoroff's formula (Brockwell and Davis, eq. 5.8.1).
In Theorems \THEOREMGLM, \THEOREMGLMB\  and \THEOREMGAP\
this result will be established
for a more general family of processes.
Since $g_Z(n)$ is the generalized variance of the process
$Z_t/{\sigma_A}$, it will be referred to as the standardized
generalized variance.
Without loss of generality we will let $\sigma_A = 1$.

{\smcaps Theorem \THEOREMFDWN .\/}
{\it Let $Z_t$ denote fractionally differenced white noise
with parameter $d \in (-{1\over2}, {1\over2})$ and $d\ne 0$.
Then $g_Z(n) \rightarrow \infty$.\/}

{\smcaps Proof.\/}
As in McLeod (1978),
$g_Z(n)=\prod_{k=0}^{n-1} \sigma_k^2$,
where $\sigma_k^2$ denotes the variance of the error in the linear
predictor of $Z_{k+1}$ using $Z_k, \ldots, Z_1$.
From the Durbin-Levinson recursion,
$$\sigma_k^2 = \cases{
  \gamma_Z(0) \qquad\qquad &  $k=0$, \cr
  \sigma_{k-1}^2 (1-\phi_{k,k}^2)\qquad\qquad  &  $k>0$. \cr}$$
where $\phi_{k,k}$ denotes the partial autocorrelation function
at lag $k$.
For the special case $p=q=0$ in \eqref{FARMA},
Hosking (1981) showed that
$\phi_{k,k} = d /(k-d)$
and
$\gamma_Z(0) =  (-2 d)! /(-d)!^2$.
Using the Durbin-Levison recursion,
$$\sigma_k^2 =
{k! (k - 2 d)! \over (k - d)!^2}
.$$
Applying the Stirling approximation to $\log(t!)$ for large $t$,
$\log (t!) \sim (t+{1\over2}) \log(t) - t + {1\over2} \log(2 \pi)$,
yields
$\log \sigma_k^2 \sim a(k)$,
where
$$a(k) = (k+{1\over 2}) \log {k (k-2 d)\over (k-d)^2}
                     + 2 d \log {k-d \over {k - 2d}}.$$
Since $\sigma_k^2$, is a monotone decreasing sequence and
for $d\ne 0$, $\sigma_k^2 > 1$,
it follows that
$\log(\sigma_k^2) $
is a positive monotone decreasing sequence.
By Stirling's approximation
$\log(\sigma_k^2)/a(k) \to 1 $ as $k \to \infty$.
So for large $k$, $a(k)$ must be a monotone decreasing sequence of
positive terms.
Expanding $a(k)$ and simplifying
$$\eqalign{a(k) &=
 (k+{1\over 2}) \log(1- {2 d \over k})
 + 2 (k+{1\over 2}) \log
    \{1+{d\over k} + \left( {d \over k} \right) ^2 + \ldots \}
 + 2 d \log (1 + { d \over k-2 d}) \cr
\noalign{\vskip6pt}
&=
 {d^2 \over k} + O({1\over k^2}),
}
$$
where the expansion $\log(1+x) = x + x^2/2 + x^3/3 + \ldots$, $|x|<1$
as been used.
Hence,
$$ k a(k) \to d^2 {\rm as\/}\  k \to \infty \eqdef{FARMAAK}$$
and by the Theorem given by Knopp (1951, \S 80, p.124),
$\sum a(k)$ diverges for $d \ne 0$.
So for $d \ne 0$, $\sum \log(\sigma_k^2)$ diverges and
consequently so does $g_Z(n)$. \BX

Eq. \eqref{FARMAAK} shows that $\sigma^2_k = 1 + O(k^{-1})$
which implies $\sigma^2_k$ decays very slowly.
The divergence of $g_Z(n)$ can be slow. See Table I.

\centerline{\smcaps Table I.}
\centerline{\smcaps Generalized variance, $g_Z(n)$, for $n=10^k,\ k=0,1,\dots,7$}
\centerline{\smcaps of fractionally differenced white noise, $Z_t$,
 with parameter $d$.}
\smallskip
$$\vbox{\halign{\strut
\hfill#\quad    
&\hfill#\quad   
&\hfill#\quad   
&\hfill#\quad   
&\hfill#\quad   
&\hfill#\quad   
&\hfill#\quad   
&\hfill#\quad   
&\hfill#\quad   
\cr             
\noalign{\hrule}
\noalign{\smallskip}
\noalign{\hrule}
\noalign{\smallskip}
$d\hfill$&$k=0$ &$k=1$&$k=2$&$k=3$&$k=4$&$k=5$&$k=6$&$k=7$  \cr
\noalign{\smallskip}
\noalign{\hrule}
\noalign{\smallskip}
$-0.4$
&1.1831&1.6225&2.3318&3.3685&4.8688&7.0375&10.1725&14.7059
\cr
$-0.1$
&1.0145&1.0366&1.0607&1.0854&1.1107&1.1365&1.1630&1.1901
\cr
$0.1$
&1.0195&1.0434&1.0678&1.0927&1.1181&1.1442&1.1708&1.1990
\cr
$0.4$
&2.0701&3.1588&4.5923&6.6417&9.6009&13.8775&20.0591&28.9951
\cr
\noalign{\smallskip}
               }}$$

\bigskip
\section{2. HYPERBOLIC DECAY TIME SERIES}

The stationary, purely nondeterministic time series, $Z_t$,
is said to be a hyperbolic decay time series with decay parameter
$\alpha$, $\alpha \in (0,2),\  \alpha \ne 1$, if for large $k$
$$
\gamma_Z(k) \sim c_{\gamma}\, k^{-\alpha},
\eqdef{Hyperbolic}
$$
where
$c_{\gamma}>0$ for $\alpha \in (0,1)$
and
$c_{\gamma}<0$ for $\alpha \in (1,2)$.
When $\alpha \in (1,2)$ the time series is said to be antipersistent.
As shown in the next theorem, antipersistent time series have
a spectral density function which decays rapidly to zero near the origin.
The term antipersistent was coined by Mandelbrot (1983)
for FGN processes with Hurst parameter, $0 <H< 1/2$.
Hyperbolic decay time series include both
FGN time series with parameter $H = 1-\alpha/2,\ H \in (0,1),\ H \ne 1/2$
and FARMA time series with parameter
$d = 1/2 - \alpha/2,\ d \in (-1/2,1/2),\ d \ne 0$.

{\smcaps Theorem \THEOREMSDAP .\/}
{\it
The spectral density function of hyperbolic decay
time series satisfies \eqref{SDFOrigin}.
}

{\smcaps Proof.\/}
Beran (1994) established this result when $\alpha \in (0,1)$ as was noted
above in eq. \eqref{SDFOrigin}.
However the Theorem of Zygmund (1968, \S V.2) used by
Beran (1994, Theorem 2.1)
does not apply to the case where $\alpha \in (1,2)$.

Let $Y_t$ have the spectral density,
$f_Y(\lambda) = c_{f} \lambda^{\alpha-1},
\ \ \alpha \in (1,2)
$.
$$\eqalign{\gamma_Y(k)
&=
2 \int\limits_0^\pi c_{f} \lambda^{\alpha-1} \cos(\lambda k) d\lambda \cr
&=
2 c_{f} k^{-\alpha} \int\limits_0^{k \pi} u^{\alpha-1} \cos(u) du, \cr
}$$
Using {\it Mathematica\/},
$$\int\limits_0^{\infty} u^{\alpha-1} \cos(u) du =
{{{\sqrt{\pi }}\,{\Gamma}({{\alpha}\over 2})}\over
    {{{({1\over 4})}^{{\alpha-1\over 2}}}\,{\rm \Gamma}({{1-\alpha}\over 2})}}$$
and so
$\gamma_Y(k) \sim c_{\gamma}k^{-\alpha},$
where
$c_{\gamma} =  2 c_{f}
{{{\sqrt{\pi }}\,{\rm \Gamma}({{\alpha}\over 2})}/
 \{
  {{({1\over 4})}^{{\alpha-1\over 2}}}\,{\rm \Gamma}({{1-\alpha}\over 2})
 \}
} < 0$.

Assume $f_Z(\lambda)$ satisfies eq. \eqref{SDFOrigin} and
we will derive \eqref{Hyperbolic}.
Since $f_Z(\lambda) / (c_f \lambda^{\alpha-1}) \to 1$ as $\lambda \to 0$,
there exists $\lambda_0$ such that for all $\lambda<\lambda_0$,
$c_f \lambda^{\alpha-1} < 1$ and
$|f_Z(\lambda) / (c_f \lambda^{\alpha-1}) - 1| < \epsilon/(2 \pi)$.
Hence for all
$\lambda<\lambda_0$,\
$|f_Z(\lambda)-f_Y(\lambda)|<\epsilon/(2 \pi)$.
Consider the systematically sampled series,
$Z_{t,\ell} = Z_{t \ell}$ for $\ell \ge 1$.
Then $Z_{t,\ell}$ has spectral density function,
$f_Z(\lambda/\ell)$.
Let $L = \pi / \lambda_0$. Then
$|f_{Z}\,(\lambda/\ell) - f_Y(\lambda)|<\epsilon/(2 \pi)$
for $\lambda \in (0,\pi)$ provided that $\ell>L$.
Hence for any $\ell>L$,
$$\eqalign{|\gamma_Z(kl)-\gamma_Y(k)|
&< 2 \int\limits_0^\pi |\cos(\lambda k)|
     |f_Z(\lambda/\ell)-f_Y(\lambda)| d\lambda \cr
&< 2 \int\limits_0^\pi
     |f_Z(\lambda/\ell)-f_Y(\lambda)| d\lambda \cr
&< \epsilon. \cr}$$
This shows \eqref{SDFOrigin} implies
\eqref{Hyperbolic}.
Since the spectral density uniquely defines the autocovariance
function, the theorem follows. \BX

Hyperbolic decay time series are self-similar: aggregated
series are hyperbolic with the same parameters as the original.

{\smcaps Theorem \THEOREMSS .\/}
{\it
Let $Z_t$ satisfy eq. \eqref{Hyperbolic} then so does $Y_t$,
where $Y_t = \sum_{j=1}^m Z_{(t-1)m+j}/m$ and $m$ is any value.
}

{\smcaps Proof.\/}
For large $\ell$,
$$\eqalign{\gamma_Y(\ell)
&=
m^{-2} \cov(\sum_{h = 1}^m Z_{(t-1) m + h},
     \sum_{k = 1}^m Z_{(t-1+\ell) m + k }) \cr
&\sim
m^{-2} \sum_{h = 1}^m \sum_{k = 1}^m c_{\gamma} (k+ m \ell-h)^{-\alpha} \cr
&\sim
m^{-2} \sum_{h = 1}^m \sum_{k = 1}^m c^{\prime}_{\gamma} \ell^{-\alpha}
(1 + {(k-h)\over m \ell})^{-\alpha} \cr
&\sim c^{\prime}_{\gamma} \ell^{-\alpha},\cr
}
$$
where $c^{\prime}_{\gamma} = m^{-\alpha} c_{\gamma}$. \BX

\bigskip
\section{3. DUALITY}

Duality has provided insights into linear time series models
(Finch, 1960;
Pierce, 1970;
Cleveland, 1972;
Box and Jenkins, 1976;
Shaman, 1976;
McLeod, 1977, 1984).
In general, the dual of the stationary invertible linear
process $Z_t = \psi(B) A_t$ is defined to be
$\psi(B) \ddot Z_t = A_t$, where
$\psi(B) = 1+\psi_1 B+\psi_2 B^2 + \ldots$
and $B$ is the backshift operator on $t$.
Equivalently, if $Z_t$ has spectral density $f_Z(\lambda)$ then
the dual has spectral density proportional to $1/f_Z(\lambda)$
with the constant of proportionality determined by the
innovation variance.
Thus in the case of a FARMA($p,q$) with parameter $d$ the dual is
a FARMA($q,p$) with parameter $-d$.
The next theorem generalizes this to the hyperbolic case.

{\smcaps Theorem \THEOREMDUAL .\/}
{\it
The dual of a hyperbolic decay time series with decay parameter $\alpha$
is another hyperbolic decay series with parameter decay parameter
$2-\alpha$.
}
\medskip
{\smcaps Proof.\/}
The spectral density near zero of the dual of a hyperbolic decay
time series with parameter $\alpha$ is
$1/(c_{f} \lambda^{\alpha-1}) = c_f^{-1} \lambda^{(2-\alpha)-1}$
which implies a hyperbolic process with parameter $2-\alpha$. \BX

{\smcaps Theorem \THEOREMHYPER .\/}
{\it The time series $Z_t$ satisfying \eqref{Hyperbolic}
may be written in random shock form as
$Z_t = A_t + \sum \psi_{\ell} A_{t-\ell} $
where
$\psi_{\ell} \sim   c_{\psi} \ell^{-(1+\alpha)/2}$
and $c_{\psi} > 0$ for $\alpha \in (0,1)$
and $c_{\psi} < 0$ for $\alpha \in (1,2)$
and in inverted form as
$Z_t = A_t + \sum \pi_{\ell} Z_{t-\ell}$
where
$\pi_{\ell} \sim    c_{\pi} \ell^{-(3-\alpha)/2}$
and $c_{\pi} > 0$ for $\alpha \in (0,1)$
and $c_{\pi} < 0$ for $\alpha \in (1,2)$
\/}
\medskip
{\smcaps Proof.\/}
The case $\alpha \in (0,1)$ was established in Theorem \THEOREMLM.
When $\alpha \in (1,2)$ the random shock coefficients are given
by
$$\eqalign{
\psi_{\ell}
&\sim
-c_{2-\pi} \ell^{-\{3-(2-\alpha)\}/2} \cr
&\sim
c_{\psi} \ell^{-(1+\alpha)/2}, \cr
}$$
where
$c_{\psi} = -c_{2-\pi}$.
Similarly for the inverted form. \BX

\bigskip
\section{4. GENERALIZED VARIANCE}
For ARMA processes, $Z_t$,
$\lim g_Z(n)$ is finite
and has been evaluated by Finch (1960) and McLeod (1977).
McLeod (1977, eq. 2) showed $g_Z(n) = m_Z + O(r^n)$,
where $r \in (0,1)$.
The evaluation of this limit uses the Theorem of
Grenander and Szeg\" o (1984, \S 5.5) which only applies
to the case where
the spectral density, $f_Z(\lambda), \lambda \in [0,2\pi)$
satisfies the Lipschitz condition
$|f_Z^{\prime}(\lambda_1)- f_Z^{\prime}(\lambda_2)|
  <K |\lambda_1-\lambda_2|^\zeta,$
for some $K>0$ and $0<\zeta<1$.
Since when $\alpha \in (0,1)$, $f_Z^{\prime}(\lambda)$ is unbounded,
this condition is not satisfied.

{\smcaps Lemma 1.\/}
{\it Let $X_t$ and $Y_t$ be any independent stationary processes
with positive  innovation variance and let $Z_t = X_t+Y_t$.
Then $G_Z(n) > G_X(n)$}

{\smcaps Proof.\/}
This follows directly from the fact that the one-step predictor
error variance of $Z_t$ can not be less than that
of $X_t$. \BX

{\smcaps Theorem \THEOREMGLM .\/}
{\it Let $Z_t$ denote a strongly persistent time process defined in
eq. \eqref{SDFOrigin}.
Then $g_Z(n) \rightarrow \infty$.\/}
\medskip

{\smcaps Proof.\/}
Since $Z_t = \sum \psi_k A_{t-k}$,
where $A_t$ is white noise with unit variance,
we can find a $q$ such that the process $Y_t$, where
$$Y_t = \sum_{k=q+1}^\infty \psi_k A_{t-k},$$
has all autocovariances nonnegative and
satisfying eq. \eqref{LongMemory}.
By using the comparison test for a harmonic series,
it must be possible to find an $N$ such that for $n>N$, the
covariance matrix $\Gamma_Y(n)$ has every row-sum greater
than $\Xi$, for any $\Xi>0$.
It then follows from Frobenius Theorem (Minc and Marcus, 1964,
p.152) that the largest eigenvalue of $\Gamma_Y(n)$ tends
to $\infty$ as $n \to \infty$.
Assume now that $\inf f_Y(\lambda) = m$ where $m>0$
and let $m_n$ denote the smallest eigenvalue of $\Gamma_Y(n)$
and let $\zeta_n$ denote the corresponding eigenvector.
Then
$$\eqalign{m_n
    & = m_n \zeta_n^\prime \zeta_n \cr
    & = \zeta_n^\prime \Gamma_Y(n) \zeta_n \cr
    & = \int\limits_{-\pi}^\pi \sum\limits_h \sum\limits_\ell
         \zeta_{n,h} \zeta_{n,\ell} e^{-i \lambda (h-\ell)}
         f(\lambda) d \lambda \cr
    & \ge 2 \pi m. \cr}$$
So $m_n \ge 2 \pi m$ and hence $g_Y(n) \to \infty$
as $n \to \infty$. By Lemma 1, $g_Z(n) \to \infty$ also.

For the more general case where $m=0$,
consider a process with spectral density function
$f(\lambda)+\epsilon$, where $\epsilon >0$.
Let $g_\epsilon(n)$ denote the standardized covariance determinant of
$n$ successive observations of this process.
So $g_\epsilon(n) \rightarrow \infty$ as $n \rightarrow \infty$
for every $\epsilon > 0$.
The autocovariance function corresponding to $f(\lambda)+\epsilon$ is
$$\gamma_\epsilon(k) = \cases{
  \gamma_Z(0) + 2 \pi \epsilon \qquad \qquad &$k = 0$,\cr
  \gamma_Z(k) \qquad \qquad &$k \ne 0$. \cr}
$$
By continuity of the autocovariance function with respect to $\epsilon$,
$\lim g_\epsilon(n) \to g_Z(n)$ as $\epsilon \to 0$. 
Let $\Xi > 0$ be chosen as large as we please
and let $\delta > 0$.
Then for any $\epsilon>0$ there exists an $N(\epsilon)$
such that for all $n \ge N(\epsilon)$,
$g_\epsilon(n) > \Xi + \delta.$
By continuity, there exists an $\epsilon_0$ such that
$g_Z(N(\epsilon_0)) > g_{\epsilon_0}(N(\epsilon_0)) - \delta.$
Hence $g_Z(N(\epsilon_0)) > \Xi$.
Since $g_Z(n+1) = g_Z(n) \sigma^2_n$,where
$\sigma^2_n > 1$ is the variance of the error of the linear predictor
of $Z_{n+1}$ given $Z_n, \ldots, Z_1$ we see that $g_Z(n)$
is nondecreasing.
It follows that $g_Z(n) > \Xi$ for all $n>N(\epsilon_0)$. \BX

Using a Theorem of Grenander and Szeg\" o (1984) this result is
easily generalized to any stationary time series, $Z_t$, for
which $\sum \gamma_Z(k) = \infty$.

{\smcaps Theorem \THEOREMGLMB .\/}
{\it Let $Z_t$ denote a time series for which
$f_Z(\lambda) \to \infty$ as $\lambda \to 0$.
Then $g_Z(n) \rightarrow \infty$.\/}
\medskip
{\smcaps Proof.\/}
From eq. (10) of Grenander and Szeg\" o (1984, \S 5.2),
as $n \to \infty$,
the largest eigenvalue of $\sigma_a^{-2} \Gamma_Z(n)$ approaches
$\sup f_Z(\lambda) = \infty$ while the smallest eigenvalue approaches
$2 \pi m$, where $m = \inf\ f(\lambda)$.
Note that Grenander and Szeg\" o's eq. (10) of \S 5.2,
applies directly to unbounded spectral densities
as is pointed by Grenander and Szeg\" o in the sentence
immediately following eq. (10), \S 5.2.
If it is assumed that $m>0$,
then the largest eigenvalue tends to infinity
and the smallest one is bounded by $2 \pi m$ as $n \rightarrow \infty$.
Hence, $g_Z(n) \rightarrow \infty$ for this special case.
The more general case where $m=0$ is handled as in Theorem \THEOREMGLM.\ \BX

In the case of ARMA models, the asymptotic covariance
determinant of the dual and primal are equal (Finch, 1960).
Since the hyperbolic decay time series are approximated
by high order AR and MA models, it might be expected
that this property holds for hyperbolic series too.
Theorem \THEOREMGAP\ which uses Lemma 2 proves that this
is the case.

{\smcaps Lemma 2.\/}
{\it
Let
$X_t = A_t + \sum_1^\infty \psi_\ell A_{t-\ell}$.
Let
$X_t(q) = A_t + \sum_1^q \psi_\ell A_{t-\ell}$,
and let $g_q(n)$ denote its standardized covariance determinant.
Then for any $\ell > 0$, $g_{q+\ell}(n) \ge g_q(n)$.
\/}

{\smcaps Proof.\/}
This follows directly from the fact that the one-step predictor
error variance of $X_t(q+\ell)$ can not be less than that
of $X_t(q)$. \BX

{\smcaps Theorem \THEOREMGAP .\/}
{\it For hyperbolic decay antipersistent time series, $Z_t$,
$g_Z(n) \rightarrow \infty$. \/}

{\smcaps Proof.\/}
Since the dual of the antipersistent time series $Z_t$
with parameter \break $2-\alpha,\  \alpha \in (0,1)$
is a strongly persistent time series $\ddot Z_t$ with parameter $\alpha$,
$\ddot Z_t$ may be represented in inverted form,
$ \ddot Z_t = A_t+\sum \pi_k \ddot Z_{t-k},$
where    $A_t$ is white noise and for large $k$,
$\pi_k \sim c_{\pi} k^{-(3 -\alpha)/2}.$
So the antipersistent time series $Z_t$ can be written,
$ Z_t = A_t - \sum \pi_k A_{t-k}.$
Let $\ddot g_L(n)$ and $g_L(n)$ denote the covariance
determinant of $n$ successive observations in the
AR($L$) and MA($L$) approximation to $\ddot Z_t$ and $Z_t$
$$ \ddot Z_t(L) = A_t+\sum_{k=1}^L \pi_k \ddot Z_{t-k}(L)$$
and
$$ Z_t(L) = A_t - \sum_{k=1}^L \pi_k A_{t-k}.$$

By Theorem \THEOREMGLM , for any $\Xi >0$ and $\delta>0$ there exists
an $N_1$ such that for $n>N_1$, $g_{\ddot Z}(n) > \Xi + \delta$.
Since $\ddot g_k(n) \to \ddot g_Z(n)$ as $k \to \infty$
there exists a $K_1(n)$ such that
$\ddot g_k(n) > \ddot g_Z(n) - \delta > \Xi$ for $k>K_1(n)$.
From McLeod (1977),
$\ddot g_k(n) = \ddot g_k(k)$ for $n \ge k$.
Hence for any $n>N_1$,
$\ddot g_k(m) > \ddot g_Z(n) -\delta > \Xi$
for $k>K_1(n)$ and $m \ge k$.
So $\ddot g_k(m) \to \infty$ as $k \to \infty$ and $m \ge k$.

Hence there exists $K_2$ such that $\ddot g_k(n) > \Xi + \delta$
for $k > K_2$ and $n \ge k$.
For any $k$, $g_k(n) = \ddot g_k(n) + O(r^n)$,
where $0<r<1$ (McLeod, 1977).
Let $k>K_2$.
Then there exists an $N_2(k)$ such that for all $n>N_2(k)$,
$g_k(n) > \ddot g_k(n) - \delta > \Xi$.
So $g_k(n) \to \infty$ as $k \to \infty$ and $n \ge k$.

For any $n$, $g_k(n) \to g_Z(n)$ as $k \to \infty$.
So for any $n$ there exists a $K_3(n)$ such that
$g_Z(n) > g_k(n) - \delta$ for all $k > K_3(n)$.
We have already established that there exists a $K_4$ such that
$g_k(n) > \Xi + \delta$ for $k>K_4$ and $n \ge k$.
Holding $n$ fixed for the moment, let $h>k$.
By Lemma 2, $g_h(n) \ge g_k(n)$.
By continuity since $h>K_4$, $g_Z(n)>g_h(n) - \delta$.
Since $g_h(n)  > \Xi + \delta  $ it follows that $g_Z(n) > \Xi$.
This establishes that $g_Z(n) \to \infty$ as $n \to \infty$. \BX

\bigskip
\section{5. CONCLUDING REMARKS}

Theorems \THEOREMGLM\ and \THEOREMGAP\ show that
hyperbolic decay time series, even antipersistent ones,
exhibit a type of long-range dependence.
The asymptotic standardized generalized variance is infinite.
This implies that the variance of the one-step linear predictor
based on the
last $k$ observations decays very slowly as compared with the
ARMA case where the decay to the innovation variance occurs
exponentially fast.
Theorem \THEOREMGLMB\ shows that this is a more general notion
of long-range dependence than the customary one.

Yakowitz and Heyde (1997) show that nonlinear Markov processes
can also exhibit strongly persistent hyperbolic decay in the
autocorrelation function.
Hence a better term for long-memory time series might be
strongly persistent hyperbolic decay series.
It is then clear that the long-range dependent aspect is merely
a characterization of the autocorrelation structure.

\bigskip

\section{REFERENCES}

\bib {\smcaps Beran, J.\/} (1994).
{\it Statistics for Long Memory Processes\/}.
London: Chapman and Hall.

\bib {\smcaps Beran, J.\/} (1995).
Maximum likelihood estimation of the differencing parameter for invertible
short and long memory autoregressive integrated moving average models.
{\it Journal of the Royal Statistical Society\/} B 57, 659--672.

\bib {\smcaps Box, G.E.P. and Jenkins, G.M.\/} (1976).
{\it Time Series Analysis: Forecasting and Control}.
(2nd edn). San Francisco: Holden-Day.

\bib {\smcaps Brockwell, P.J. and Davis, R.A.\/} (1991).
{\it Time Series: Theory and Methods\/}.
(2nd edn.). New York: Springer-Verlag.

\bib {\smcaps Cleveland, W.S.\/} (1972).
The inverse autocorrelations of a time series and
their applications.
{\it Technommetrics\/} 14,   277--293.


\bib {\smcaps Diebold, F.X. \& Rudebusch, G.D.\/} (1989).
Long memory and persistence in aggregate output.
{\it Journal of Monetary Economics\/}, 24, 189--209.

\bib {\smcaps Finch, P.D.\/} (1960).
On the covariance determinants of autoregressive and
moving average models.
{\it Biometrika\/} 47, 194--196.

\bib {\smcaps Graham, R.L., Knuth, D.E. and Patashnik, O.} (1989).
{\it Concrete Mathematics\/}.
Reading: Addison-Wesley.

\bib {\smcaps Granger, C.W.J. and Joyeux, R.\/} (1980).
An introduction to long-range time series models and fractional
differencing.
{\it Journal of Time Series Analysis\/} 1, 15--30.

\bib {\smcaps Grenander, U. and Szeg\" o, G.\/} (1984).
{\it Toeplitz Forms and Their Applications\/},
(2nd edn.). New York: Chelsea.

\bib {\smcaps Hipel, K.W. \& McLeod, A.I.\/} (1994).
{\it Time Series Modelling of Water Resources and Environmental Systems}.
Elsevier: Amsterdam.

\bib {\smcaps Hosking, J.R.M. \/} (1981).
Fractional differencing.
{\it Biometrika\/} 68, 165--176.

\bib {\smcaps Knopp, K. \/} (1951).
{\it Theory and Application of Infinite Series\/},
(2nd edn.). New York: Hafner.

\bib {\smcaps Mandelbrot, B.B.M.\/} (1983).
{\it The Fractal Geometry of Nature\/}.
San Francisco: Freeman.

\bib {\smcaps McLeod, A.I.\/} (1977).
Improved Box-Jenkins estimators.
{\it Biometrika\/} 64, 531--534.

\bib {\smcaps McLeod, A.I.\/} (1984).
Duality and other properties of multiplicative seasonal
autoregressive-moving average models.
{\it Biometrika\/} 71, 207--211.

\bib {\smcaps Minc, H. and Marcus, M} (1964).
{\it A Survey of Matrix Theory and Matrix Inequalities\/}.
Boston: Prindle, Weber and Schmidt.

\bib {\smcaps Pierce, D.A.\/} (1970).
A duality between autoregressive and moving average
processes concerning their parameter estimates.
{\it Annals of Statistics\/} 41, 722--726.

\bib {\smcaps Shaman, P.\/} (1976).
Approximations for stationary covariance matrices and their
inverses with application to ARMA models.
{\it Annals of Statistics\/} 4, 292--301.

\bib {\smcaps Yakowitz, S.J. and Heyde, C.C.\/} (1997).
Long-range dependency effects with implications for forecasting
and queueing inference.
Unpublished Manuscript.

\bib {\smcaps Zygmund, A.\/} (1968).
{\it Trigonometric Series\/}.
London: Cambridge University Press.

\end{document}